\documentclass[12pt]{report}
\usepackage{amssymb,amsmath,makeidx,verbatim}
\newcommand{\be}{\begin{enumerate}}
	\newcommand{\ee}{\end{enumerate}}
\newcommand{\bq}{\begin{eqnarray*}}
	\newcommand{\eq}{\end{eqnarray*}}

\begin{document}
	\newcommand{\disp}{\displaystyle}
	\thispagestyle{empty}
	\begin{center}
		\textsc{On the Lie algebra of a Malcev algebra\\}
		\ \\
		\textsc{Olufemi O. Oyadare}\\
		\ \\
		Department of Mathematics,\\
		Obafemi Awolowo University,\\
		Ile-Ife, $220005,$ NIGERIA.\\
		\text{E-mail: \textit{femi\_oya@yahoo.com}}\\
		
	\end{center}
	\begin{quote}
		{\bf Abstract.} {\it This paper develops the structure theory of a Malcev algebra via the consideration of its most important and largest Lie (sub-) algebra. We introduce the notion of a Lie algebra which uniquely corresponds to a Malcev algebra and use this correspondence to derive some basic properties of some types of ideals in the Malcev algebra. We then prove the exact nature of the root-space decomposition of a Malcev algebra. }
	\end{quote}
	\ \\
	
	{\bf $\bf{\S 1.}\;$ Introduction.}\\
	
	One of the most understood and well-developed class of algebras in the whole of mathematics is the class of {\it Lie algebras.} A Lie algebra is generally defined as an algebra $\mathfrak{L}$ (over $\mathbb{R}$ or $\mathbb{C}$) whose product, $[\cdot,\cdot],$ is anti-commutative (i.e., $[x,x]=0$) and in which the Jacobian, $$[[x,y],z]+[[y,z],x]+[[z,x],y]=:J(x,y,z),$$ vanishes, for all $x,y,z\in\mathfrak{L}.$ Numerous articles, monographs and textbooks by notable authors have been published detailing different aspects of a Lie algebra, their inter-relationships and applications to the life sciences. We shall refer to Varadarajan $[10.]$ for orientation. Perhaps, the most singular reason for the extreme utility of a Lie algebra, over and above other algebras, stems out of the fact that every Lie algebra corresponds to a Lie group. This affords results to be exchanged between a Lie algebra (with its {\it infinitesimal} theory) and its corresponding Lie group (with its {\it topological} and {\it manifold} structures), with one enriching the other in a uniquely determined way. Indeed no other class of algebras has had the privilege of the class of Lie algebras.
	
	Since the {\it Malcev algebra} (also called the {\it Moufang-Malcev algebra}) was introduced (See Malcev $[5.],$ Schafer $[9.]$ and Paal $[6.]$) as a tangent space of an analytic {\it Moufang loop} and since it was shown that every Lie algebra is a Malcev algebra, mathematicians have been seeking the Lie-algebra type of approach to the understanding of a Malcev algebra. It is believed that a Malcev algebra is a higher version of a Lie algebra in which $J=J(x,y,z)$ is not necessarily zero. Aspects of the structure theory of a Malcev algebra had been discussed by A. A. Sagle in his University of California thesis $[8.].$ It had also been found that it is possible to have a form of the Poincare-Birkhoff-Witt Theorem (known for a Lie algebra) for a Malcev algebra, as well as its (nonassociative) universal enveloping algebra (See $[2.],$ $[7.]$ and $[9.]$). All these partial results ($[3.]$) were however not enough to give the class of Malcev algebras the esteemed status occupied by one of its sub-classes, the (sub-) class of Lie algebras.
	
	The present paper gives a novel approach to the study of a Malcev algebra via the consideration of a distinguished Lie (sub-) algeba closely attached to it. The paper is organised as follows: \S 2 considers the structure theory of a Malcev algebra while \S develops its basic ideal theory. The root-space decomposition of a Malcev algebra is the subject of \S4.
	\ \\
	
	{\bf $\bf{\S 2.}$ Structure theory of a Malcev algebra}\\
	
	We shall start by giving a precise definition to a Malcev algebra and then give a rapid presentation of aspects of its structure.\\
	\indent {\bf 2.1 Definition.} A {\it Malcev algebra} $A$ over a field $k$ is any nonassociative algebra over $k$ endowed with a distributive multiplication $[\cdot,\cdot]$ satisfying the relations $$[a,a]=0$$ and $$[J(a,b,c),a]=J(a,b,[a,c]),$$ for all $a,b,c\in A.\;\Box$
	
	\indent It should be clear from the above definition that a Lie algebra is in particular a Malcev algebra: Indeed, $[J(a,b,c),a]=[0,a]=0=J(a,b,[a,c]).$ One could then say that a Malcev algebra is a Lie algebra whose Jacobian may not necessarily vanish but satisfies $[J(a,b,c),a]=J(a,b,[a,c]).$ This relationship assures us of many examples of a Malcev algebra while examples of a Malcev non-Lie algebra could be found in the literature, Sagle $[8.].$
	
	\indent The result that got the subject started is the following.\\
	\indent {\bf 2.2 Lemma.} ($[8.],$ p. $435$) Let $A$ denote a Malcev algebra and let $J(A,A,A)$ denote its linear subspace spanned by all elements of the form $J(a,b,c)$ (for $a,b,c\in A$). Then $J(A,A,A)$ is an ideal of $A\;\Box.$
	
	\indent We shall often simply write $J$ for the ideal $J(A,A,A).$ It should be noted that for a Lie algebra, $A,$ the ideal $J(A,A,A)$ is $\{0\}.$ It then follows that the residue-class algebra $A/J$ is well defined. In the special case of a Lie algebra $A$ we have that $$A/J=A/\{0\}\simeq A.$$ It happens that $J(A,A,A)$ is the most important ideal for a Malcev non-Lie algebra. Indeed, we have the following.\\
	\indent {\bf 2.3 Lemma.} ($[8.],$ p. $435$) Let $A$ denote a Malcev algebra of characteristic not $2$ or $3.$ Then $J(A,A,A)$ is the smallest ideal of $A$ for which $A/J$ is a Lie algebra$.\;\Box$
	
	\indent This in particular shows that $J$ is uniquely determined as an ideal of $A.$ Another distinguished subset of a Malcev algebra is defined in the following.\\
	\indent {\bf 2.4 Definition.} Let $A$ denote a Malcev algebra. The {\it $J-$nucleus} of $A$ is given as $$N=\{x\in A: J(x,A,A)=0\}.\;\Box$$
	
	\indent This means that $N$ is the maximal subset of $A$ for which $J(N,A,A)=0$ and may also be shown to an ideal of $A.$ Hence $N$ is uniquely determined as an ideal of $A.$ These two uniqely determined ideals $J$ and $N$ of the Malcev algebra $A$ annihilates each other.\\
	\indent {\bf 2.5 Lemma.} ($[8.],$ p. $440$) If $A$ is a Malcev algebra, then $NJ=0.\;\Box$
	
	\indent The third most important ideal of a Malcev algebra is also given below.\\
	\indent {\bf 2.6 Definition.} Let $A$ denote a Malcev algebra and define $\triangle(x,y)$ by the requirement $$[z,\triangle(x,y)]=J(z,x,y),$$ for all $x,y,z\in A.$
	
	\indent Elementary properties of $\triangle$ may be found in Sagle$[8.]$ while the crucial result concerning $\triangle$ is the following.\\
	\indent {\bf 2.7 Theorem.} ($[8.],$ p. $432$) Let $A$ denote a Malcev algebra of characteristic not $2$ or $3$ and let $\triangle(A,A)$ denote the linear span of all the $\triangle(x,y)$ (for all $x,y\in A$). Then $\triangle(A,A)$ is a Lie algebra under commutation$.\;\Box$\\
	
	\indent The last result shows that $\triangle(A,A)$ is not realized as a Lie sub-algebra of $A,$ since it is only a Lie algebra with respect to the commutation $$(x,y)=[x,y]-[y,x],$$ $x,y\in A.$ Hence, it is bound to exhibit properties not known to $A.$ Moreover, by Proposition $5.8$ of Sagle $[8.],$ the Lie algebra, $\triangle(A,A)$ is not uniquely determined from $A.$
	
	\indent Our approach in this paper is to employ the well-known fine structure of the Lie (sub-)algebra $N$ of a Malcev algebra $A$ to get across to all of $A$ via the residue-class Lie algebra $A/J.$ The first major result along this line is the following.\\
	\indent {\bf 2.8 Theorem.} ($[8.],$ p. $442$) Let $A$ denote a Malcev algebra such that $\triangle(A,A)$ is completely reducible in $A,$ then $$A=N\bigoplus J,$$ where $N$ is completely reduced by $\triangle(A,A)$ and $J$ is a semi-simple ideal of $A.\;\Box$\\
	
	\indent The hypotheses of Theorem $1.8$ hold if $A$ is semi-simple. This theorem gives a very useful relationship between the uniquely determined ideals $N$ and $J$ (of $A$) in the light of the residue-class Lie algebra $A/J.$ This is our point of departure in this paper.\\

	{\bf $\bf{\S 3.}$ Ideal theory of a Malcev algebra}\\
	
	We use the truth of Theorem $2.8$ to define the Lie algebra of a Malcev algebra.\\
	\indent {\bf 3.1 Definition.} Let $A$ denote a Malcev algebra. The Lie algebra corresponding to or of $A$ is defined as the residue-class (Lie algebra) $A/J.\;\Box$\\
	
	\indent {\bf 3.2 Lemma.} Every Malcev algebra has a unique Lie algebra corresponding to it.\\
	\indent {\bf Proof.} Let $A$ denote any Malcev algebra, then the $J-$nucleus $N$ is non-empty and it is the maximal subset of $A$ for which $J(N,A,A)=0.$ The direct sum $A=N\bigoplus J$ of Theorem $2.8$ thus implies that $A/J\simeq N.$ These respectively imply the existence and uniqueness of the Lie algebra $A/J.\;\Box$\\
	
	\indent In the case where $A$ is a semi-simple Malcev algebra, its corresponding Lie algebra $A/J$ is equally semi-simple. We shall denote by $\phi$ the canonical epimorphism, $A:\rightarrow A/J$ and shall often make use of the following basic result of {\it commutative algebra.}\\
	
	\indent {\bf 3.3 Lemma.} ($[1.],$ p. $2$) Let $A$ denote any Malcev algebra. There is a one-to-one order preserving correspondence between the ideals $\mathfrak{b}$ of $A$ which contain $J$ and the ideals $\overline{\mathfrak{b}}$ of the Lie algebra $A/J,$ given by $\mathfrak{b}=\phi^{-1}(\overline{\mathfrak{b}}),$ in which prime ideals correspond to prime ideals$.\;\Box$\\
	
	\indent It therefore follows from the last lemma that, even though there are examples to show that it is not always true that the product of ideals of a MAlcev algebra $A$ is an ideal of $A$ (which is in contrast to the fine property of ideals in a Lie algebra), this failure could now be given a proof and could also be salvaged for ideals of $A$ contain $J$ as contained as follows.
	
	\indent {\bf 3.4 Proposition.} If $\mathfrak{a}$ and $\mathfrak{b}$ are any two ideals of a Malcev algebra, $A,$ both of which contain $J,$ then $\mathfrak{a}\mathfrak{b}$ is an ideal of $A.$\\
	\indent {\bf Proof.} The one-to-one correspondence of the last lemma implies that $\mathfrak{a}=\phi^{-1}(\overline{\mathfrak{a}})$ and $\mathfrak{b}=\phi^{-1}(\overline{\mathfrak{b}}),$ for some ideals $\overline{\mathfrak{a}}$ and $\overline{\mathfrak{b}}$ in the Lie algebra $A/J.$ Thus $$\mathfrak{a}\mathfrak{b}=\phi^{-1}(\overline{\mathfrak{a}})\phi^{-1}(\overline{\mathfrak{b}})
	=\phi^{-1}(\overline{\mathfrak{a}}\overline{\mathfrak{b}}),$$ where it is known that $\overline{\mathfrak{a}}\overline{\mathfrak{b}}$ is an ideal in $A/J\simeq N.$ The crucial Lemma $3.3$ thus implies that $\mathfrak{a}\mathfrak{b}$ is an ideal in $A$ containing $J.\;\Box$\\
	
	\indent Our approach in this paper will further show that the true correspondence between a Malcev algebra and a Lie algebra is via the process of its residue-class Lie algebra $A/J,$ which will make the beautiful results of a Lie algebra available to the residue-class (Lie) algebra and then to the Malcev algebra. This could be gleaned through the following commutative diagram: $$A\rightarrow^{\phi} A/J\rightarrow^{i} N,$$ where $\phi$ is the canonical map of Lemma $3.3$ and $i$ is the isomorphism of Lemma $3.2.$ It then means that we have the projection $$(i \circ \phi):A\rightarrow N.$$ Clearly, $i(N+J)=N.$
	
	\indent The above result, for ideals of $A$ containing $J,$ may now be extended to include any subset of the Lie algebra $A/J\simeq N,$ which would buttress our point of view that it is most profitable to approach the study of a Malcev algebra via its projection map $(i \circ \phi):A\rightarrow N,$ into the Lie algebra $N.$\\
	\indent {\bf 3.5 Lemma.} Let $\mathfrak{h}$ denote a subalgebra (respectively, an ideal) of $N,$ then $H:=(i \circ \phi)(\mathfrak{h})$ is a subalgebra (respectively, an ideal) of $A.$\\
	\indent {\bf Proof.} The surjectivity of $\phi$ implies that for any $n\in N$ there is $a\in A$ such that $a=(i \circ \phi)(n).\;\Box$
	
	\indent Now that we have a direct access between $A$ and its Lie algebra, $N,$ it would be instructive if the extra condition (included in Proposition $3.4$) of an ideal containing $J$ can be removed. This now takes us to the consideration of a special class of ideals of $N$ meant for this purpose.
	
	\indent {\bf 3.6 Definition.} Let $A$ denote a Malcev algebra. A subset $X$ of the Lie algebra $N$ is called an {\it $i-$ideal of $N$} whenever $$i(n)X\subseteq X,$$ for every $n\in A/J.\;\Box$
	
	\indent Lemma $3.5$ (together with Lemma $3.3$) suggests that the requirement of an $i-$ideal $X$ of $N$ is equivalent to having $(i\circ \phi)(a)X\subseteq X,$ for every $a\in A.$ The following result shows the central importance of $i-$ideal of the Lie algebra $N$ in getting access to and in the understanding of {\it all} the ideals of a Malcev algebra.\\
	\indent {\bf 3.7 Theorem.} Let $A$ denote any Malcev algebra. Then $\mathfrak{a}$ is an ideal of $A$ iff $(i \circ \phi)(\mathfrak{a})$ is an $i-$ideal of $N.$\\
	\indent {\bf Proof.} Let $a\in A/J$ and let $y\in (i \circ \phi)(\mathfrak{a}),$ then $y=(i \circ \phi)(x),$ for some $x\in \mathfrak{a},$ an ideal of $A.$ It is clear (by Lemma $3.5$) that $X:=(i \circ \phi)(\mathfrak{a})$ is an additive abelian subgroup of $N$ and that $(i \circ \phi)(ax)\subseteq X.$ Note also that $$(i \circ \phi)(ax)=i(\phi(ax))=i(a\phi(x))=i(a)(i \circ \phi)(x)=i(a)y\subseteq i(a)X$$ and that $$i(a)(i \circ \phi)(x)=(i \circ \phi)(ax)\subseteq X.$$ Hence, $i(a)X\subseteq X,$ for every $a\in N.$
	
	\indent Conversely, let $\overline{\mathfrak{a}}$ denote an $i-$ideal of $N.$ That is, let $i(a)\overline{\mathfrak{a}}\subseteq \overline{\mathfrak{a}},$ for every $a\in A/J.$ Consider now a subset $\mathfrak{a}$ of $A$ as $$\mathfrak{a}:=\{(\phi^{-1}\circ i^{-1})(y):y \in \overline{\mathfrak{a}}\}.$$ Clearly, $(i\circ\phi)(\mathfrak{a})=\overline{\mathfrak{a}}.$ Now let $r\in A$ and let $x\in\mathfrak{a},$ then $x=(\phi^{-1}\circ i^{-1})(y),$ for some $y\in \overline{\mathfrak{a}}.$ Considering the element $rx,$ we have that $(i \circ \phi)(rx)=i(\phi(rx))=i(r\phi(x))=i(r)(i \circ \phi)(x)\in i(r)(i \circ \phi)(\mathfrak{a})\in i(r)\overline{\mathfrak{a}}\subseteq \overline{\mathfrak{a}},$ for all $r\in A.$ In particular, for all $r\in N.$ Hence, $(i \circ \phi)(rx)=i(\phi(rx))\in \overline{\mathfrak{a}},$ for all $r\in N.$ Whence, $rx\in (\phi^{-1}\circ i^{-1})(\overline{\mathfrak{a}})=\mathfrak{a}.\;\Box$
	
	\indent The last Theorem sets up a one-to-one correspondence between (all of) the ideals of a Malcev algebra $A$ and the $i-$ideals of its Lie algebra $N.$ It equally confirms that the Lie algebra $N$ is replete with different types of ideals, while the type which is directly concerned with a Malcev algebra is the $i-$ideals of $N.$ We now give a general proof (in contrast to the usual counter-examples) of a major difference on the behaviour of ideals in a Malcev algebra and in a Lie algebra. The next two results show that this negative behaviour of ideals in a Malcev algebra is due to the properties of $i-$ideals in $N.$\\
	\indent {\bf 3.8 Proposition.} The product of two $i-$ideals of $N$ may not always be an $i-$ideal of $N.$\\
	\indent {\bf Proof.} Every $i-$ideal of $N$ is in particular an ideal in $A.$ Hence, if $\overline{\mathfrak{a}}$ and $\overline{\mathfrak{b}}$ are $i-$ideals of $N$ they also satisfy $\overline{\mathfrak{a}}\overline{\mathfrak{b}}\subseteq\overline{\mathfrak{a}}$ and $\overline{\mathfrak{a}}\overline{\mathfrak{b}}\subseteq\overline{\mathfrak{b}}.$ Thus $$i(a)(\overline{\mathfrak{a}}\overline{\mathfrak{b}})\subseteq i(a)\overline{\mathfrak{a}}\subseteq\overline{\mathfrak{a}}$$ and $$i(a)(\overline{\mathfrak{a}}\overline{\mathfrak{b}})\subseteq i(a)\overline{\mathfrak{b}}\subseteq\overline{\mathfrak{b}},$$ for any $a\in A/J.$ Therefore, $$i(a)(\overline{\mathfrak{a}}\overline{\mathfrak{b}})\subseteq\overline{\mathfrak{a}}\bigcap\overline{\mathfrak{b}}\nsubseteq\overline{\mathfrak{a}}\overline{\mathfrak{b}}.\;\Box$$
	
	\indent The following may therefore be seen as the first known general result on the product of ideals of a Malcev algebra $A.$\\
	\indent {\bf 3.9 Theorem.} Product of ideals of a Malcev algebra $A$ is not always an ideal of $A.$\\
	\indent {\bf Proof.} We simply combine Theorem $3.7$ with Proposition $3.8.\;\Box$\\
	
	\indent Since there is a one-to-one correspondence between ideals of $A$ with the $i-$ideals of the Lie algebra $N,$ it follows that among the $i-$ideals $\overline{\mathfrak{a}}$ of $N$ are those (that form a special class) that corresponds to the ideals $\phi^{-1}(\overline{\mathfrak{a}})$ of $A$ that contain $J$ and therefore satisfy Proposition $3.4.$ Basic properties of operations on ideals suggest that we seek these (special ones) among the co-prime $i-$ideals of $N$ in order to complete the proof of Proposition $3.8.$\\
	\indent {\bf 3.10 Lemma.} The product of co-prime $i-$ideals of $N$ is an $i-$ideal of $N.$\\
	\indent {\bf Proof.} For every $a\in N$ and for co-prime $i-$ideals $\overline{\mathfrak{a}}$ and $\overline{\mathfrak{b}}$ of $N,$ $$i(a)(\overline{\mathfrak{a}}\overline{\mathfrak{b}})\subseteq\overline{\mathfrak{a}}\bigcap\overline{\mathfrak{b}}\subseteq\overline{\mathfrak{a}}\overline{\mathfrak{b}}.\;\Box$$
	
	\indent It is very interesting to confirm if the co-prime $i-$ideals of $N$ correspond to ideals of $A$ containing $J.$\\
	
	\textbf{\S 4. Root-space decomposition of a Malcev algebra}
	
	We now consider how the root-space decomposition of the Lie algebra $N$ enriches $A.$\\
	\indent {\bf 4.1 Definition.} Let $A$ denote a Malcev algebra and let $N$ be its corresponding Lie algebra. The {\it adjoint action} of $A$ on $N$ is the map $$ad_{A\times N}:A\times N\rightarrow N$$ $$:(a,n)\rightarrowtail ad_{A\times N}(a,n)$$ $$=:ad_{A\times N}(a)(n)$$ defined as $$ad_{A\times N}(a)(n)=i(ad_{A/J}\phi(a)(i^{-1}(n))).$$
	
	\indent The map $ad_{A/J}:A/J\rightarrow A/J$ is the adjoint map on the Lie algebra $A/J,$ given as $ad_{A/J}(x)(y)=[x,y].$ The following is immediate.\\
	\indent {\bf 4.2 Lemma.} The adjoint action of $A$ on $N$ is the adjoint map on $N.$\\
	\indent {\bf Proof.} $$ad_{A\times N}(a)(n)=i(ad_{A/J}\phi(a)(i^{-1}(n)))=i([\phi(a),i^{-1}(n)])=[(i\circ\phi)(a),n].\;\Box$$
	\indent One of the key results of a Lie algebra is the following.\\
	\indent {\bf 4.3 Lemma.} ($[4.],$ p. $132$) Let $A$ denote a finite-dimensional Malcev algebra over $\mathbb{C}$ and let $\mathfrak{h}$ be a nilpotent Lie subalgebra of the Lie algebra $N,$ then the generalised weight spaces $N_{\alpha}$ of $N$ relative to $ad_{N}\mathfrak{h}$ satisfy\\
	\indent $(a)$ $N=\bigoplus N_{\alpha},$ where $N_{\alpha}$ is the set $$\{X\in N: (adH-\alpha(H)\cdot 1)^{n}X=0,\;\forall\;H\in\mathfrak{h},\mbox{for some integer}\;n=n(H,X)\geq0\};$$
	\indent $(b)$ $\mathfrak{h}\subseteq N_{o};$\\
	\indent $(c)$ $[N_{\alpha},N_{\beta}]\subset N_{\alpha+\beta}$
	(with the understanding that $N_{\alpha+\beta}$ is $0$ if $\alpha+\beta$ is not a generalised weight).\\
	\indent We then use these properties of $N_{\alpha}$ in $N$ to develop an equivalent theory for generalised weight spaces of a finite-dimensional Malcev algebra, $A.$\\
	\indent {\bf 4.4 Theorem.} Let $A$ denote a finite-dimensional Malcev algebra over $\mathbb{C}$ and let $N_{\alpha}$ be the generalised weight spaces $N_{\alpha}$ of $N$ relative to $ad_{N}\mathfrak{h}$ (given in Lemma $4.3$ above). Choose a subset $H$ of $A$ such that $\mathfrak{h}:=(i\circ\phi)(H)$ is a nilpotent Lie subalgebra of the Lie algebra $N.$ Then the generalised weight spaces $A_{\alpha}$ of $A$ relative to $ad_{A\times N}H$ satisfy\\
	\indent $(a)$ $A=\bigoplus A_{\alpha},$ where $A_{\alpha}$ is given as the set $$\{a\in A: (ad_{A\times N}h-\alpha(h)\cdot 1)^{n}a=0,\;\forall\;h\in\mathfrak{h},\mbox{for some integer}\;n=n(h,a)\geq0\};$$
	\indent $(b)$ $H\subseteq A_{o};$\\
	\indent $(c)$ $[A_{\alpha},A_{\beta}]\subset N_{\alpha+\beta}$
	(with the understanding that $N_{\alpha+\beta}$ is $0$ if $\alpha+\beta$ is not a generalised weight).\\
	\indent {\bf Proof.} We set $A_{\alpha}=(i\circ\phi)(N_{\alpha}).$ Hence, $$[A_{\alpha},A_{\beta}]=[(i\circ\phi)(N_{\alpha}),(i\circ\phi)(N_{\beta})]
	=i[\phi(N_{\alpha}),\phi(N_{\beta})]$$ $$=i[N_{\alpha}+J,N_{\beta}+J] =i([N_{\alpha},N_{\beta}]+J)\subseteq i(N_{\alpha+\beta}+J)=N_{\alpha+\beta},$$ since $N_{\alpha+\beta}\subseteq N,$ in which $i(N+J)=N.\;\Box$\\
	\indent There is still a lot of the results of Lie algebras (and Lie groups) that may be proven for a Malcev non-Lie algebra.\\
	
	{\bf \S 5. References.}
	\begin{description}
		\item [{[1.]}] Atiyah, M. F. and Macdonald, I. G., \textit{Introduction to commutative algebra,} Addison-Wesley Series in Mathematics, Addison-Wesley, Reading. $1969.$
		\item [{[2.]}] Bremner, M. R, Hentzel, I. R., Peresi, L. A., Tvalavadze, M. V., Usefi, H., {\it Enveloping algebras of Malcev algebras,} {\it Comment. Math. Univ. Carolin.} {\bf 51}, $2$ ($2010$), p. $157-174.$
		\item [{[3]}] Harrathi, F., Mabrouk, S., Ncib, O., Silvestrov, S. {\it Malcev Yang-Baxter equation, weighted $\mathfrak{O}-$operators on Malcev algebras and post-Malcev algebras,} arXiv:$2206.03629v1$ [math.RA] ($2022$)
		\item [{[4.]}] Knapp, A. W. {\it Lie Groups Beyond an Introduction,} Progress in Mathematics, {\bf 140} ($2002$)
		\item [{[5.]}] Malcev, A. I., Analytic loops, {\it Matem. Sb.} {\bf 78}, ($1955$), p. $569-578.$
		\item [{[6.]}] Paal, E., Note on analytic Moufang loops, {\it Comment. Math. Univ. Carolin.} {\bf 45}, $2$ ($2004$), p. $349-354.$
		\item [{[7.]}] Perez-Izquiredo, J. M., Shestakov, I. P., An envelope for Malcev algebras, {\it J. Algebra} $\bf 272$ ($2004$), p. $379-393.$
		\item [{[8.]}] Sagle, A. A., Malcev algebras, Trans. Amer. Math. Soc. {\bf 101} ($1961$), p. $426-458.$
		\item [{[9.]}] Schafer, R. D., {\it An introduction to nonassociative algebras,} Academic Press, New York, ($1966$).
		\item [{[10.]}] Varadarajan, V. S., {\it Lie groups, Lie algebras and their representations,} Springer-Verlag, New York, ($1984$).
	\end{description}
	
\end{document}